\documentclass[12pt,a4paper]{article}
\usepackage[top=4.2cm, bottom=2.8cm, left=3cm, right=3cm]{geometry}
\usepackage{mathrsfs}

\usepackage[francais, english]{babel}
\usepackage[utf8]{inputenc}
\usepackage[T1]{fontenc}
\usepackage{amsfonts, amssymb, amsthm}
\usepackage{dsfont}
\usepackage{thmtools}
\usepackage{mathtools}
\usepackage{amsmath}

\usepackage{multirow}
\usepackage{graphicx}
\usepackage{caption, subcaption}
\usepackage{float}
\usepackage{enumitem}

\usepackage{pict2e}

\usepackage{mdframed}

\usepackage[all]{xy}

\usepackage{tikz}

\usepackage{multicol} %Environnement Multicolones

\mathtoolsset{showonlyrefs} %Numerotation Equations ssi Citee 

\usepackage{relsize,exscale}

\usepackage{hyperref}
\usepackage{url}
 {\begin{list}{}%
         {\setlength{\leftmargin}{#1}}%
         \item[]%
 }
 {\end{list}}

\allowdisplaybreaks %autorisation de casser les environements align pour changer de page

\newcounter{TheoremA}
\newtheorem{ThA}[TheoremA]{Theorem}
\newtheorem{PropA}[TheoremA]{Proposition}

\newtheorem*{Cor*}{Corollary} %Corollaire non numéroté pour l'intro

\newtheorem{Th}{Theorem}[section] %Le compteur recommence à chaque section
\newtheorem{Prop}[Th]{Proposition}

\theoremstyle{definition}
\newtheorem{Def}[Th]{Definition}

\theoremstyle{remark}
\newtheorem{Remark}[Th]{Remark}

\newcounter{fig}

\newcommand{\Lie}{\mathcal{L}}
\newcommand{\R}{\mathbb{R}}
\newcommand{\C}{\mathbb{C}}
\newcommand{\scal}[2]{\left\langle#1,#2\right\rangle}

\DeclareMathOperator{\Vol}{Vol}
\newcommand{\Cinf}{\mathcal{C}^{\infty}}
\newcommand{\Hspace}{\mathcal{H}}

\newcommand{\poisson}[2]{\left\lbrace#1,#2\right\rbrace}
\newcommand{\Sph}{{\mathbb{S}}}

\title{Isometries of the Space of Sasaki Potentials}
\author{Thomas Franzinetti}

\begin{document}
\maketitle
\abstract{
Given any two Kähler manifolds $X_1$ and $X_2$, L. Lempert recently proved that if their spaces of Kähler potentials are isometric with respect to the Mabuchi metric, then $X_1$ and $X_2$ must be diffeomorphic. We prove that this is no longer the case for Sasaki manifolds. Then, considering \emph{regular} Sasaki manifolds $M_1$ and $M_2$, we prove that if the spaces of potentials are isometric, then $M_1$ and $M_2$ must have, among others, the same universal covering space. Finally, getting rid of the regularity assumption on $M_1$ and $M_2$, we investigate the consequences of the existence of \emph{affine} Mabuchi isometries: this leads to a family of Sasaki isospectral structures.
}

\section{Introduction}

Let $(X_i, \omega_i)$, $i=1,2$ be Kähler manifolds of complex dimension $n$; one can consider the space of Kähler potentials:
$$
\mathcal{H}(X_i, \omega_i) = \left\lbrace u \in \Cinf(X) \; \vert \; \omega_{i,u} := \omega_i + dd^c u >0 \right\rbrace.
$$ 
These infinite dimensional spaces are furnished with a Riemannian metric, the so called Mabuchi metric \cite{mabuchi_symplectic_1986}:
$$
\scal{\phi_1}{\phi_2}^{\mathcal{H}(X_i, \omega_i)}_u := \frac{1}{\Vol(X_i, \omega_i)} \int_{X_i} \phi_1 \phi_2 \; \omega_{i,u}^n \; , \quad \phi_1,\phi_2 \in \Cinf(X_i), u\in\mathcal{H}(X_i, \omega_i),
$$
where $\Vol(X_i, \omega_i) := \int_{X_i}\omega_i^n$. In this setting, we have the following classification result of isometries with respect to the Mabuchi metric.

\begin{ThA}[{\cite[Theorem 1.1]{darvas2019isometries} \& \cite[Theorem 1.1]{lempert2017isometries}}]
\label{Thm:INTROIsometriesKahler}
Let $u\in \mathcal{H}(X_1, \omega_1)$ and suppose $F : \mathcal{H}(X_1, \omega_1) \rightarrow \mathcal{H}(X_2, \omega_2)$ is a global isometry. Then there are a diffeomorphism $G:X_2\rightarrow X_1$, and a real number $b$ such that: $G^*\left(\omega_{1,u}\right) = \pm \omega_{2,F(u)}$ and
$$
F_*\vert_u (\phi) = G^*\phi - b\int_{X_1}\phi \; \omega_{1,u}^n, \; \quad \phi \in T_u\mathcal{H}(X_1,\omega_1) \simeq \Cinf(X_1).
$$
The possible values for $b$ are $0$ or $2/\Vol(X_1, \omega_1)$.
\end{ThA}

Note that \cite[Theorem 1.1]{darvas2019isometries}, which is a consequence of \cite[Theorem 1.1]{lempert2017isometries}, is actually proved in the case $(X_1,\omega_1) = (X_2,\omega_2)$. However, the proof can be adapted word to word in the more general setting of arbitrary two $(X_1, \omega_1)$ and $(X_2,\omega_2)$. In addition, the map $G$ can be proved to be a biholomorphism or an anti-biholomorphism \cite[Lemma 3.8]{darvas2019isometries}.

It is worth mentioning using the Monge-Ampère flip \cite[Section 2.3]{darvas2019isometries}, one can see that every $G$ and $b$ as above are realized by isometries.

This theorem allows us to identify the underlying Kähler manifolds provided their spaces of potentials are isometric. Thus Theorem \ref{Thm:INTROIsometriesKahler} turns \cite[Theorem 1.1]{darvas2019isometries} in a sort of Banach-Stone theorem for Mabuchi isometries in the sense that some geometrical properties of a Kähler manfiold can be read at the level of its space of potentials: the existence of an isometry $\mathcal{H}(X_1, \omega_1) \rightarrow \mathcal{H}(X_2, \omega_2)$ for the Mabuchi metrics automatically implies that the underlying manifolds $X_1$ and $X_2$ must be diffeomorphic and the diffeomorphism $X_2 \rightarrow X_1$ pulls back the Kähler class $\left[\omega_1 \right]$ to $\left[\omega_2 \right]$. Besides, this Theorem shows that the differential $F_\star$ depends on the point where it is calculated only in the case $b\neq 0$.

\medskip

A Sasaki structure $\mathcal{S} = (\xi, \eta, \Phi, g)$ on a smooth manifold $M$ it the data of a contact form $\eta$ with its Reeb vector field $\xi$, an almost complex structure $\Phi$ (i.e. endomorphisms on each tangent spaces satisfying $\Phi^2 = -1 + \xi\otimes \eta$) and a Riemannian metric $g$. All these tensors must satisfy compatibility conditions:
$$
\Phi\circ\Phi = -\mathds{1}_{TM} + \xi \otimes \eta \; ; \quad g\circ(\Phi\otimes \mathds{1}_{TM}) = \frac{1}{2}d\eta \; ; \quad g\circ(\Phi\otimes\Phi) = g - \eta\otimes\eta.
$$
We refer to the excellent monograph \cite{boyer_sasakian_2007} concerning Sasaki geometry.

Given a Sasaki structure $\mathcal{S} = (\xi, \eta, \Phi, g)$ on a smooth $(2n+1)-$dimensional manifold $M$, one can define an analogue of the Mabuchi metric on the space of Sasaki potentials \cite{guan_geodesic_nodate, guan_regularity_2009}:
$$
\mathcal{H}(M, \mathcal{S}) = \left\lbrace u \in \Cinf_B(M) \; \vert \; d\eta_{u} := d\eta + dd^c u >0 \right\rbrace,
$$
where $\Cinf_B(M)$ is the space of basic functions identified with the tangent space of $\mathcal{H}(M, \mathcal{S})$ at any point. The Riemannian metric on $\Hspace(M,\mathcal{S})$ is then:
$$
\scal{\phi_1}{\phi_2}_u^{\mathcal{H}(M,\mathcal{S})} = \frac{1}{\Vol(M, \mathcal{S})}\int_{M} \phi_1 \phi_2 \; \eta \wedge d\eta_{u}^n \; , \quad \phi_1,\phi_2 \in \Cinf_B(M), u\in\mathcal{H}(M, \mathcal{S}),
$$
where $\Vol(M, \mathcal{S}) := \int_{M} \eta \wedge \left( d\eta\right)^n$.

Sasaki geometry is often considered as the odd dimensional analogue of Kähler geometry. To that extend, one can expect a result similar to Theorem \ref{Thm:INTROIsometriesKahler}. However, we have:

\begin{PropA}[Proposition \ref{Prop:OurExample}]
\label{PropINTRO:OurExample}
There exist Sasaki structures on the 3-sphere $\Sph^3$ and on quotients $\Sph^3 / \mathbb{Z}_k$ for any $k\in \mathbb{N}^*$ such that the corresponding spaces of potentials are isometric.

In particular, in the Sasaki setting, the existence of an isometry between the spaces of potentials doesn't imply that the underlying manifolds are diffeomorphic.
\end{PropA}

Recall that a Sasaki structure $(M,\mathcal{S})$ is said to be \emph{regular} if the action generated by the Reeb flow is free. If the action is only locally free, the Sasaki structure is said to be \emph{quasi-regular} otherwise, it is said to be \emph{irregular}.

A regular (resp. quasi-regular) Sasaki manifold is, among others, the total space of a circle bundle over $(X,\omega)$, a Kähler manifold (resp. orbifold -- see \cite[Section 4.3]{boyer_sasakian_2007}). The Kähler class $\left[\omega\right]$ represents the characteristic class of the circle bundle. In particular, $\left[ \omega \right]$ is an integral class and $X$ is said to be a \emph{Hodge manifold} (resp. \emph{Hodge orbifold}). We refer to \cite[Theorem 7.1.3 \& Theorem 7.5.1]{boyer_sasakian_2007} for more details.

\smallskip

As illustrated by Proposition \ref{PropINTRO:OurExample}, the existence of an isometry between spaces of potentials is less restrictive in the Sasaki setting. However we still have rigidity results. The first one concerns \emph{regular} Sasaki manifolds. 

\begin{ThA}[Theorem \ref{Th:RegularSasakiIsometryH}]
\label{Th:INTRORegularSasakiIsometryH}
Given any two regular Sasaki manifolds $(M_i, \mathcal{S}_i)$ over $(X_i,\omega_i),\; i=1,2$ such that there is an isometry:
$$
F : \mathcal{H}(M_1, \mathcal{S}_1) \longrightarrow \mathcal{H}(M_2, \mathcal{S}_2)
$$
Then there must be a diffeomorphism $G : X_2 \rightarrow X_1$ such that $G^*\left[\omega_1\right] = \pm \left[ \omega_2 \right]$. In particular, $M_1$ and $M_2$ must have the same universal covering space. Moreover, if $l_i$ is the period of $(M_i,\mathcal{S}_i)$ (that is the period of the flow of $\xi_i$), $i=1,2$, then:
\begin{equation}
%\label{eq:RatioRationnel}
\frac{l_2}{l_1} \in \mathbb{Q}.
\end{equation}
\end{ThA}

%%%%%%%%%%%%%%%%%%%%%%%%%
\bigskip

To any $u\in \Hspace(M,\mathcal{S})$, one can associate a new Sasaki structure $\mathcal{S}_u$ whose contact form is $\eta_u = \eta + d^c u$, the Reeb vector field is unchanged. Its almost complex structure is $\Phi_u := \Phi - \xi \otimes(d^c u \circ \Phi)$ and the Riemannian metric is $g_u := \eta_u \otimes \eta_u + \frac{1}{2} d\eta \circ (\mathop{1}_{TM} \otimes \Phi_u)$. 

Now, we can relax the regularity assumption and get a weaker result concerning \emph{affine isometries}. By \emph{affine isometry}, we mean that the differential $F_\star\vert_u$ is independent on $u$. As illustrated by Theorem \ref{Thm:INTROIsometriesKahler}, in the Kähler setting, affine isometries are the one with $b=0$.

\begin{ThA}[Theorem \ref{Th:AffineIsometries}]
\label{Th:INTRO_AffineIsometries}
Let  $F : \mathcal{H}(M_1, \mathcal{S}_1) \longrightarrow \mathcal{H}(M_2, \mathcal{S}_2)$ be an affine isometry. Then, for any $u\in\Hspace(M_1,\mathcal{S}_1)$, the basic Laplacians of $\left(M_1, \mathcal{S}_{1,u}\right)$ and the one of $\left(M_2, \mathcal{S}_{2,F(u)}\right)$ are isospectral.
\end{ThA}

We observe that this Theorem does not need any regularity assumption on the considered Sasaki manifolds, meaning that we even allow us to work with irregular Sasaki manifolds.

The basic Laplacian, defined in Section \ref{sec:AffineIsometriesSasaki}, can be simply seen as the restriction of the Riemannian Laplacian to the space of basic functions. 

In particular, if $(M_i, \mathcal{S}_i)$ is a quasi-regular Sasaki manifold over the Kähler orbifold $(X_i,\omega_i),\; i=1,2$, then the underlying orbifolds must have an infinite dimensional family of pairwise isospectral metrics: for each potential $u$, the Kähler metric associated to $\mathcal{S}_{i,u}$ must be isospectral to the Kähler metric associated to $\mathcal{S}_{2,F(u)}$. Interpreting Theorem \ref{Th:INTRO_AffineIsometries} in terms of Kähler orbifolds, one can state a weak form of Theorem \ref{Thm:INTROIsometriesKahler}:

\begin{Cor*}
Given any two Hodge orbifolds whose associated circle bundle is a Sasaki manifold, assume that there is an affine isometry between spaces of potentials, then there must be many pairwise isospectral metrics on the orbifolds.
\end{Cor*}

The circle bundle associated to a Hodge orbifold is a Sasaki manifold if the local uniformizing groups \cite[Definition 4.1.1]{boyer_sasakian_2007} are subgroups of $\mathbb{S}^1$ \cite[Theorem 7.5.2]{boyer_sasakian_2007}.

In the case of smooth Kähler manifolds, Theorem \ref{Thm:INTROIsometriesKahler} gives, in the same way as above, many pairwise \emph{isometric} metrics.

%%%%%%%%%%%%%%%%%%%%%%%%%
\vspace{1cm}

We organize this paper as follows: Section \ref{sec:RegularCase} is devoted to regular Sasaki manifolds, we give an example proving Proposition \ref{PropINTRO:OurExample}, then after a short preliminary concerning circle bundles, we give a detailed proof of Theorem \ref{Th:INTRORegularSasakiIsometryH}. Then, in Section \ref{sec:AffineIsometriesSasaki}, we prove Theorem \ref{Th:INTRO_AffineIsometries}.

\section*{Acknowledgments}
I am grateful for the constant support of Eleonora Di Nezza and Gilles Courtois without whom this work would have never emerged. I also thank Tam{\'a}s Darvas, Eveline Legendre and L{\'a}szl{\'o} Lempert for listening to me and answering to my questions. Their valuable advice were really helpful to me.

\section{Isometries of the Space of Sasaki Potentials: the Regular Case}
\label{sec:RegularCase}

We will focus on the first kind of Sasaki structures: the regular ones. Let $(M, \mathcal{S}) = (M, \xi, \eta, \Phi, g)$ be a regular Sasaki manifold. Each orbit of $\xi$ has the same period $l$ \cite[Section I.2.]{BoothbyWang1958}. It is well known that the space of orbits $X$ is a Kähler manifold and one can chose the Kähler form $\omega$ such that the quotient projection $\pi : M \rightarrow X$ pulls back $\omega$ on $\frac{1}{2} d\eta$ (see for instance \cite[Section 7.5]{boyer_sasakian_2007}). In this case, the volumes verify:
$$
\Vol(M,\mathcal{S}) = \int_M \eta\wedge d\eta^n = 2^n \int_M \eta \wedge \left(\frac{d\eta}{2}\right)^n = 2^n l \Vol(X, \omega).
$$

Let us describe the link between the infinite dimensional Riemannian spaces $\mathcal{H}(M,\mathcal{S})$ and $\mathcal{H}(X,\omega)$. Given any basic function $\overline{\phi} \in \Cinf_B(M)$, there is a unique smooth function $\phi\in\Cinf(M)$ such that $2\pi^*\phi = \overline{\phi}$. It is straightforward to check that the map $2\pi^* : \mathcal{H}(X,\omega) \rightarrow \mathcal{H}(M,\mathcal{S})$ is a linear continuous diffeomorphism which is identified with its differential map $\Cinf(X) \rightarrow \Cinf_B(M)$. We claim that it is not only a diffeomorphim but also a homothetic transformation with respect to the Mabuchi metrics. Indeed, given $\overline{u} = 2\pi^* u \in \mathcal{H}(M,\mathcal{S})$ and $\overline{\phi_j} = 2\pi^* \phi_j \in T_{\overline{u}} \mathcal{H}(M,\mathcal{S}) \simeq \Cinf_B(M)$ for $j=1,2$, we have:
\begin{align*}
\Vol(M,\mathcal{S})\scal{\overline{\phi_1}}{\overline{\phi_2}}_{\overline{u}}^{\mathcal{H}(M,\mathcal{S})}
&= \int_M \overline{\phi_1}\,\overline{\phi_2} \; \eta\wedge d\eta_{\overline{u}}^n \\
&=\int_M 2\pi^*(\phi_1) \cdot 2\pi^*(\phi_2) \; \eta \wedge \left(d\eta + dd^c(2\pi^* u) \right)^n \\
&= 2^{n+2}\int_M \pi^*(\phi_1 \phi_2) \; \eta\wedge \left( \frac{d\eta}{2} + dd^c(\pi^* u) \right)^n \\
&= 2^{n+2}\int_M \eta\wedge \pi^* \left( \phi_1 \phi_2 \omega_u^n \right) \\
&= 2^{n+2} l \int_X \phi_1 \phi_2 \; \omega_u^n.\\
\end{align*}
From this, we deduce:
\begin{align}
\scal{\overline{\phi_1}}{\overline{\phi_2}}_{\overline{u}}^{\mathcal{H}(M,\mathcal{S})} 
&= \frac{2^{n+2} l}{\Vol(M,\mathcal{S})} \int_X \phi_1 \phi_2 \; \omega_u^n\\
&= \frac{2^{n+2} l}{2^n l \Vol(X,\omega)} \int_X \phi_1 \phi_2 \; \omega_u^n\\
&= 4\scal{\phi_1}{\phi_2}_u^{\mathcal{H}(X,\omega)}. \label{eq:computationLinkMabuchiSasakiKahler}
\end{align}

\bigskip

These notions are all we need to prove Proposition \ref{PropINTRO:OurExample}:

\begin{Prop}%[Proposition \ref{PropINTRO:OurExample}]
\label{Prop:OurExample}
There exist Sasaki structures on the 3-sphere $\Sph^3$ and on quotients $\Sph^3 / \mathbb{Z}_k$ for any $k\in \mathbb{N}^*$ such that the corresponding spaces of potentials are isometric. In particular, in the Sasaki setting, the existence of an isometry between the spaces of potentials doesn't imply that the underlying manifolds are diffeomorphic.
\end{Prop}

\begin{proof}
Consider the Hopf fibration over the $3$-sphere (endowed with its standard Sasaki structure $\mathcal{S}$). The circle action $\Sph^1 \curvearrowright \Sph^3 \subset \C^2$ is given by
$$
e^{i\theta}\cdot (z_1, z_2) := (e^{i\theta} z_1, e^{i\theta} z_2).
$$
For any $k\in \mathbb{N}$, the subgroup $\mathbb{Z}_k = \mathbb{Z}/k\mathbb{Z} \subset \Sph^1$ acts on $\Sph^3$ leaving the Sasaki structure invariant. Thus, we get a family of regular Sasaki manifolds $(M_k := \Sph^3 / \mathbb{Z}_k, \mathcal{S}(k))$. Fix $k\in \mathbb{N}^*$. Both $\Sph^3$ and $M_k$ come with a projection $\pi_{\Sph^3} : \Sph^3 \rightarrow \mathbb{CP}^1$ and $\pi_{M_k} : M_k \rightarrow \mathbb{CP}^1$. Using twice \eqref{eq:computationLinkMabuchiSasakiKahler}, it is clear that the composition $F_0$ of the maps:

\begin{equation}
\label{eq:IsometryExampleSphere}
\begin{array}{cccccc}
F_0: &\mathcal{H}(M_k, \mathcal{S}(k)) & \rightarrow & \mathcal{H}(\mathbb{CP}^1, \omega_{FS}) & \rightarrow & \mathcal{H}(\Sph^3, \mathcal{S}_a) \\
 & 2\pi_{M_k}^* u & \mapsto & u & \mapsto & 2\pi_{\Sph^3}^* u
\end{array}
\end{equation}
is an isometry. In this example, the underlying manifolds $\mathbb{S}^3$ and $M_k$ are not even homeomorphic: their first fundamental group are different.  However, both $\Sph^3$ and $M_k$ are circle bundles over $\mathbb{CP}^1$ furnished with the Fubini-Study metric. Moreover, they have the same universal covering space.
\end{proof}

%main Thm regular Sasaki
\bigskip
We will show that the situation above is somewhat generic. The following Theorem which is Theorem \ref{Th:INTRORegularSasakiIsometryH} of the introduction is a consequence of Theorem \ref{Thm:INTROIsometriesKahler}.

\begin{Th}%[Theorem \ref{Th:INTRORegularSasakiIsometryH}]
\label{Th:RegularSasakiIsometryH}
Given any two regular Sasaki manifolds $(M_i, \mathcal{S}_i)$ over $(X_i,\omega_i)$, $i=1,2$ such that there is an isometry:
$$
F : \mathcal{H}(M_1, \mathcal{S}_1) \longrightarrow \mathcal{H}(M_2, \mathcal{S}_2)
$$
Then there must be a diffeomorphism $G : X_2 \rightarrow X_1$ such that $G^*\left[\omega_1\right] = \pm \left[ \omega_2 \right]$. In particular, $M_1$ and $M_2$ must have the same universal covering space. Moreover, if $l_i$ is the period of $(M_i,\mathcal{S}_i)$ (that is the period of the flow of $\xi_i$), $i=1,2$, then:
\begin{equation}
\label{eq:RatioRationnel}
\frac{l_2}{l_1} \in \mathbb{Q}.
\end{equation}
\end{Th}

\begin{Remark}
It follows form the proof that the last two conditions in Theorem \ref{Th:RegularSasakiIsometryH} are consequences of the fact that $M_1$ and $M_2$ must have the same transverse space with same Kähler class. In addition, it follows from the examples of Proposition \ref{Prop:OurExample} that every ratios of periods as in \eqref{eq:RatioRationnel} are possibles.
\end{Remark}
As an immediate corollary, if $(M,\mathcal{S})$ is a regular Sasaki manifold, then for any positive real number $\lambda\neq 1$, the spaces $\mathcal{H}(M, \mathcal{S})$ and $\mathcal{H}(M,\mathcal{S}_\lambda)$ (see below for a definition of $\mathcal{S}_\lambda$) are not isometric.

%%%% Circle bundles
\bigskip
Before proving Theorem \ref{Th:INTRORegularSasakiIsometryH}, we will do a brief recap on circle bundles. Given any positive real parameter $\lambda>0$, we write $\mathcal{S}_\lambda$ for the $\lambda-$\emph{transverse homothetic deformed} Sasaki structure \cite[(7.3.10)]{boyer_sasakian_2007}. Namely, if $\mathcal{S} = (\xi, \eta, \Phi, g)$, then
\begin{equation}
\label{eq:DefTransverseHomotheticSasakian}
\mathcal{S}_\lambda := \left(\frac{1}{\lambda} \xi, \lambda\eta, \Phi, g_\lambda := \lambda g + (\lambda^2-\lambda)\eta\otimes\eta\right).
\end{equation}
If $l$ it the period of the Reeb flow of $(M,\mathcal{S})$, the Sasaki manifold $(M, \mathcal{S}_{\frac{1}{l}})$ is the total space of a $\R/\mathbb{Z}-$bundle over $X$. Such bundles are classified by their characteristic class \cite[Section 3]{kobayashi1956}. Here, we write $e(M,\mathcal{S})$ for the characteristic class of the underlying circle bundle of the Sasaki manifold $(M,\mathcal{S})$. In this specific case this characteristic class is proportional to $\omega$ \cite[Theorem 2]{BoothbyWang1958}. More precisely:
\begin{equation}
\label{eq:FirstChernClassCircleBundle}
e(M,\mathcal{S})_b = \left[ \frac{2}{l}\omega \right] \in H^{1,1}(X)\cap H^2(X,\mathbb{Z})_b.
\end{equation}
Here, the subscript $b$ is the notation of S. Kobayashi for the image of the natural map
$$
\xymatrix{
H^2(X,\mathbb{Z})\ar@{->>}[r]& H^2(X,\mathbb{Z})_b \ar@{}[r]|-*[@]{\subset}& H^2(X,\mathbb{R}).
}
$$
S. Kobayashi proved \cite[Theorem 2]{kobayashi1956} that the map sending the total space $M$ of a circle bundle to its characteristic class $e(M,\mathcal{S})_b$ is a surjective group map (the group structure is defined in \cite[p30 \& Theorem 1]{kobayashi1956}) and if $e(M,\mathcal{S})_b = 0$, then the bundle $\overline{M}$, defined in \cite[p37]{kobayashi1956}, induced from the universal covering space $\tilde{X} \rightarrow X$ must be trivial \cite[p33]{kobayashi1956}. In particular, if the space of orbits $X$ is simply connected, $M$ must be a trivial bundle. Extensive details about these definitions, the group structure and circle bundles can be found in \cite{kobayashi1956}. We will use this fact later to identify two circle bundles. Note that \eqref{eq:FirstChernClassCircleBundle} implies that for a regular Sasaki manifold $(M,\mathcal{S})$, one cannot expect $e(M,\mathcal{S})_b = 0$.

\begin{proof}[Proof of Theorem \ref{Th:INTRORegularSasakiIsometryH}]
For $i = 1,2$, let $l_i$ be the period of $\xi_i$ on $M_i$ and let $X_i$ be the space of orbits with quotient projection $\pi_i$ and Kähler class $\omega_i$. We will first find out what does the existence of an isometry F imply at the level of Kähler potentials. Consider the following diagram defining $\Psi$:
$$
\xymatrix{
	\mathcal{H}(M_1, \mathcal{S}_1) \ar[r]^{F} &\mathcal{H}(M_2,\mathcal{S}_2) \\
	\mathcal{H}(X_1, \omega_1) \ar[u]^{2\pi_1^*}\ar[r]_\Psi & \mathcal{H}(X_2, \omega_2) \ar[u]_{2\pi_2^*}
}
$$
Again, using twice \eqref{eq:computationLinkMabuchiSasakiKahler}, we immediately get that $\Psi$ is an isometry. Thanks to Theorem \ref{Thm:INTROIsometriesKahler}, we get a diffeomorphism $G: X_2 \rightarrow X_1$ such that:
\begin{equation}
\label{eq:diffeoLempert}
\left[G^*(\omega_1)\right] = \varepsilon\left[\omega_2 \right] \in H^{1,1}(X_2),\text{ where } \varepsilon = \pm 1.
\end{equation}
In particular:
$$
\frac{l_1}{2} G^*\left(e(M_1, \mathcal{S}_1)_b\right) = \varepsilon\frac{l_2}{2} e(M_2, \mathcal{S}_2)_b.
$$
Since the characteristic classes are integral, we can get \eqref{eq:RatioRationnel}:
\begin{equation}
\frac{l_1}{l_2} =: \frac{p}{q} \in \mathbb{Q}\, , \quad p,q\in \mathbb{N}.
\end{equation}
Assume first $\varepsilon = 1$. Let $G^*(M_1)$ be the circle bundle over $X_2$ induced by $G$. Since $G$ is a diffeomorphism, as circle bundles we identify $M_1 \simeq G^*(M_1)$. As the total space of a circle bundle over $X_2$, $p \, G^*(M_1) - q \, M_2$ has null characteristic class. Thus, the bundle $\overline{p\, G^*(M_1) - q\, M_2}$ induced by the universal covering space $\tilde{X_2}\rightarrow X_2$ must be trivial. According to \cite[Theorem 4]{kobayashi1956}, the map sending a bundle to the induced one by a fixed diffeomorphism is a group morphism. So $p\overline{G^*(M_1)} - q \overline{M_2}$ is trivial. Using \cite[Theorem 5]{kobayashi1956}, we get the following series of identifications as circle bundles:
$$
\overline{M_1}/\mathbb{Z}_p \simeq p\overline{M_1} \simeq q \overline{M_2} \simeq \overline{M_2}/\mathbb{Z}_q.
$$
We claim that $\overline{M_i}\rightarrow M_i$, $i=1,2$ is a covering. %cfp174 Mars2020

\begin{proof}[Proof of the claim]
Let $\pi : P \rightarrow X$ be a circle bundle over $X$. Let $p:\tilde{X} \rightarrow X$ be the universal covering space of X. Consider the induced circle bundle over $\tilde{X}$, namely:
$$
\xymatrixcolsep{5pc}
\xymatrix{
\overline{P} := \left\lbrace (x,q) \in \tilde{X}\times P \; \vert \; \pi(q) = p(x) \right\rbrace \ar[r]^-{\text{proj}_2} \ar[d]_{\text{proj}_1} & P \ar[d]^\pi\\
\tilde{X} \ar[r]^p & X
}
$$
Let $q \in P$, and choose a small neighbourhood $\mathcal{U}$ of $q$ in $P$ such that $\pi(\mathcal{U})$ is contained in a local trivialization $\mathcal{V}$ of $\tilde{X} \rightarrow X$: $p^{-1}(\mathcal{V}) = \sqcup_i \mathcal{V}_i$ and $p\vert_{\mathcal{V}_i}$ is a homeomorphism form $\mathcal{V}_i$ to $\mathcal{V}$. We have the disjoint union of open subsets of $\overline{P}$:
$$
\text{proj}_2^{-1}(\mathcal{U}) = \bigsqcup_i \left((\mathcal{V}_i \times \mathcal{U}) \cap \overline{P} \cap \text{proj}_2^{-1}(\mathcal{U}) \right).
$$
Each of these subsets maps homeomorphically to $\mathcal{U}$:
\begin{equation}
\begin{array}{ccccc}
\mathcal{U}& \rightarrow &  \left((\mathcal{V}_i \times \mathcal{U}) \cap \overline{P} \cap \text{proj}_2^{-1}(\mathcal{U}) \right)\\
r  & \mapsto &  \left( p\vert_{\mathcal{V}_i}^{-1}(\pi(r)), r \right)
\end{array}
\end{equation}
whose reciprocal map is $\text{proj}_2$.
\end{proof}

Denoting $\tilde{M_i}$ the universal covering space of $M_i$ and $\tilde{\overline{M_i}}$ the universal covering space of $\overline{M_i}$, $i=1,2$, we have the following diagram where each arrow is a covering map:
$$
\xymatrix{
 \tilde{M_1} \ar[d] & \tilde{\overline{M_1}}\ar[d] & \tilde{\overline{M_2}}\ar[d] & \tilde{M_1}\ar[d]\\
 M_1 & \overline{M_1}\ar[l]\ar[d] & \overline{M_2}\ar[r]\ar[d] & M_2 \\
 & {\overline{M_1}}/{\mathbb{Z}_p} \ar@{}[r]|-*[@]{\simeq} & {\overline{M_2}}/{\mathbb{Z}_q} &
}
$$
The central part of the diagram shows $\tilde{\overline{M_1}} \simeq \tilde{\overline{M_2}}$ while the left hand part and the right hand part respectively show $\tilde{M_1} \simeq \tilde{\overline{M_1}}$ and $\tilde{M_2} \simeq \tilde{\overline{M_2}}$. In particular, the universal covering spaces of $M_1$ and $M_2$ must be the same, proving the theorem if $\varepsilon = 1$.

In the case $\varepsilon = -1$, one can restart the argument with the conjugated Sasaki structure on $M_2$ for instance. Namely, we consider $M_2$ with the Sasaki structure $\mathcal{S}_2^c := (-\xi_2, -\eta_2, -\Phi_2, g_2)$ defined in \cite[Definition 7.5.11]{boyer_sasakian_2007}. Then, as infinite dimensional Riemannian manifolds, the space of potentials is unchanged. Indeed, under such a change, the orientation is reversed, the operator $d^c$ is changed to its opposite (we write $d^{c,\mathcal{S}_2} = -d^{c,\mathcal{S}^c_2}$ with obvious notations) and the notion of positivity for basic $(1,1)$-forms is reversed. Thus,
\begin{align*}
\phi \in \mathcal{H}(M_2, \mathcal{S}_2)
& \iff d\eta + dd^{c,\mathcal{S}_2} \phi \text{ is } \mathcal{S}_2\text{-positive } \\
& \iff - d\eta - dd^{c,\mathcal{S}_2} \phi  = - d\eta + dd^{c,\mathcal{S}^c_2} \phi \text{ is } \mathcal{S}_2 \text{-negative } \\
& \iff - d\eta + dd^{c,\mathcal{S}^c_2} \phi \text{ is } \mathcal{S}^c_2 \text{-positive } \\
& \iff \phi \in \mathcal{H}(M_2, \mathcal{S}^c_2).
\end{align*}
Since the orientation on $M_2$ is defined with the volume form $\eta_2\wedge d\eta_2^n$, it is changed by $(-1)^{n+1}$. This shows that neither the volume nor the Mabuchi metric change: $\forall \phi_1, \phi_2 \in \Cinf_B(M_2), \forall u \in \mathcal{H}(M,\mathcal{S})$,

\begin{align*}
\scal{\phi_1}{\phi_2}_u^{\mathcal{H}(M_2,\mathcal{S}^c_2)} 
&= \frac{(-1)^{n+1}}{\Vol(M_2,\mathcal{S}_2)} \int_{M_2} \phi_1 \phi_2 (-\eta)\wedge \left(-d\eta + dd^{c,\mathcal{S}^c_2}\right)^n \\
&= \frac{(-1)^{n+1}}{\Vol(M_2,\mathcal{S}_2)} \int_{M_2} \phi_1 \phi_2 (-\eta)\wedge \left(-d\eta - dd^{c,\mathcal{S}_2}\right)^n \\
&= \scal{\phi_1}{\phi_2}_u^{\mathcal{H}(M_2, \mathcal{S}_2)}.
\end{align*}
Here, all the integral are considered on $M_2$ with the orientation defined by $\mathcal{S}_2$. The isometry $F : \mathcal{H}(M_1, \mathcal{S}_2) \rightarrow \mathcal{H}(M_2, \mathcal{S}_2)$ induced by the diffeomorphism $G$ gives an isometry $\mathcal{H}(M_1, \mathcal{S}_2) \rightarrow \mathcal{H}(M_2, \mathcal{S}^c_2)$ induced by a diffeomorphism as in \eqref{eq:diffeoLempert} with $\varepsilon = 1$ since the Kähler class $\left[\omega_2\right]$ is changed to its opposite.
\end{proof}

%%%%%%%%%%%%%%%%%%%%%%%%%%%%%%%%%%%%%%%%%%%%%%%%%%%%%%%%%%%%%%%%%%%%%%%%%%%%%%%%%%%%%%%%%%%%%%%%%%%%%%%%%%%%%%%%%%%%%%%%%%%%%%%%%%%%%
%%%%%%%%%%%%%%%%%%%%%%%%%%%%%%%%%%%%%%%%%%%%%%%%%%%%%%%%%%%%%%%%%%%%%%%%%%%%%%%%%%%%%%%%%%%%%%%%%%%%%%%%%%%%%%%%%%%%%%%%%%%%%%%%%%%%%
%%%%%%%%%%%%%%%%%%%%%%%%%%%%%%%%%%%%%%%%%%%%%%%%%%%%%%%%%%%%%%%%%%%%%%%%%%%%%%%%%%%%%%%%%%%%%%%%%%%%%%%%%%%%%%%%%%%%%%%%%%%%%%%%%%%%%
%%%%%%%%%%%%%%%%%%%%%%%%%%%%%%%%%%%%%%%%%%%%%%%%%%%%%%%%%%%%%%%%%%%%%%%%%%%%%%%%%%%%%%%%%%%%%%%%%%%%%%%%%%%%%%%%%%%%%%%%%%%%%%%%%%%%%
%%%%%%%%%%%%%%%%%%%%%%%%%%%%%%%%%%%%%%%%%%%%%%%%%%%%%%%%%%%%%%%%%%%%%%%%%%%%%%%%%%%%%%%%%%%%%%%%%%%%%%%%%%%%%%%%%%%%%%%%%%%%%%%%%%%%%

\section{Affine Isometries: Regular, Quasi-Regular and Irregular case.}
\label{sec:AffineIsometriesSasaki}
A natural question would be to find a generalization Theorem \ref{Th:INTRORegularSasakiIsometryH} in the case of quasi-regular Sasaki manifolds (or even irregular ones). Recall that a quasi-regular Sasaki manifold is a circle bundle over a Kähler orbifold \cite[Theorem 7.1.3]{boyer_sasakian_2007}, thus we would like to find an analogue of \cite[Theorem 1.1]{lempert2017isometries} for orbifolds. However, even in the quasi-regular case, L. Lempert's techniques cannot be used. Indeed, L. Lempert applies Lie algebra techniques for the Poisson bracket \cite[Section 4]{lempert2017isometries} since the curvature of the Mabuchi metric can be expressed with a triple Poisson bracket \cite[Theorem 4.3]{mabuchi_symplectic_1986}. Thus, since an isometry preserves the curvature operator, it also preserves the triple bracket \cite[Lemma 4.1]{lempert2017isometries}. The key point, is to show that such isometries preserve, up to a sign, the Poisson bracket itself \cite[Lemma 4.4]{lempert2017isometries}. This is a consequence of the fact that any function of vanishing integral can be written as the sum of some Poisson brackets \cite[Lemma 4.6]{lempert2017isometries}. This is no longer true on orbifolds as we explain in what follows.

\medskip

Recall that the Poisson bracket of any two functions $\phi_1, \phi_2$ can be expressed as an exact derivative:
$$
\poisson{\phi_1}{\phi_2} = \left(\nabla^c \phi_1 \right) \cdot \phi_2,
$$
where $\nabla^c$ stands for the Hamiltonian vector field (sometimes called symplectic gradient). On a Kähler orbifold \cite[Section 4.4.4]{boyer_sasakian_2007}, $\nabla^c \phi = J\nabla \phi$ where $J$ is the complex structure and $\nabla$ is the Riemannian gradient. Now consider a singular point $p$ with local uniformizing group $\Gamma_p$ \cite[Definition 4.1.1]{boyer_sasakian_2007}. Note that by \cite[Corollary 7.1.7]{boyer_sasakian_2007}, $\Gamma_p$ must be cyclic in our setting. Recall that by the very definition \cite[Definition 4.2.9]{boyer_sasakian_2007}, a vector field must be $\Gamma_p$-invariant at $p$. In particular, if $\Gamma_p$ contains a rotation, then any vector field must vanish at $p$. As a consequence, any Poisson bracket must also vanish at $p$. Now, it is clear that there are functions of vanishing integral which are \emph{not} vanishing at $p$. In particular, such functions can not be written as the sum of Poisson brackets.

\medskip

Our idea is then to use the fact that an isometry preserves the covariant derivative.

\bigskip

Preservation of the curvature by an isometry, which is a third order condition, is the cornerstone in \cite{lempert2017isometries}. Here, we will take advantage of the fact that an isometry preserves the covariant derivative which is a second order condition. We don't make any assumptions about the regularity of the Sasaki manifolds $(M_i, \mathcal{S}_i), \; i = 1,2$. But we restrict ourselves to specific isometries.

\begin{Def}
We say that an isometry is an \emph{affine isometry} if its differential map is independent of the point where it is calculated.
\end{Def}

Note that, if $F$ is an affine map, for any $\phi_0,\phi_1 \in \Hspace$, we write $\phi_t := t\phi_1 + (1-t)\phi_0$ and we have:
\begin{align}
F(\phi_1) &= F(\phi_0) + \int_0^1 \frac{d}{dt}F(\phi_t) dt \\
		&= F(\phi_0) + \int_0^1 F_\star(\phi_1 - \phi_0)dt \\
		&= F(\phi_0) + F_\star(\phi_1-\phi_0) \label{eq:AffineMap}.
\end{align}

In the Kähler case, affine isometries correspond to the isometries with $b=0$ in Theorem \ref{Thm:INTROIsometriesKahler}. To that extend, it seems reasonable to restrict ourselves to affine isometries.

The end of this section is devoted to the proof of Theorem \ref{Th:INTRORegularSasakiIsometryH} of the introduction:

\begin{Th}
\label{Th:AffineIsometries}
Let  $F : \mathcal{H}(M_1, \mathcal{S}_1) \longrightarrow \mathcal{H}(M_2, \mathcal{S}_2)$ be an affine isometry. Then, for any $u\in\Hspace(M_1,\mathcal{S}_1)$, the basic Laplacians $\Delta_{u,B}$ on $M_1$ and $\Delta_{F(u),B}$ on $M_2$ are isospectral.
%Moreover, $F_\star (1) = 1$.

In particular, if $M_1$ and $M_2$ are quasi-regular Sasaki manifolds, then the underlying orbifolds must have many isospectral metrics.
\end{Th}

As it is well known for Riemannian manifolds, there are examples of non-isometric isospectral orbifolds. Such examples are provided for instance in \cite{weilandt2012isospectral} where isospectral families of pairwise non-isometric metrics are constructed on weighted projective spaces.

\bigskip

Let us assume now that there is an affine isometry 
$$
F : \mathcal{H}(M_1, \mathcal{S}_1) \longrightarrow \mathcal{H}(M_2, \mathcal{S}_2)
$$
between spaces of potentials over Sasaki manifolds.

Since $F$ is an isometry, it preserves not only the curvature (i.e. the triple bracket) but also the covariant derivative \cite[Lemma 5.1]{lempert2017riemannian}. Let us recall the expression of the covariant derivative \cite[Definition 3.2]{franzinetti2020curvature}: given a Sasaki manifold $(M,\mathcal{S})$, to each $u\in \Hspace(M,\mathcal{S})$, one can associate a new Sasaki structure $\mathcal{S}_u$ as in the introduction. We denote by $\nabla$ the gradient with respect to $g_u$. Now, consider $s\mapsto u(s)$ a path in $\Hspace(M,\mathcal{S})$ and $s \mapsto Y(s)\in \Cinf_B(M)$ a vector field tangent to $\Hspace(M,\mathcal{S})$ along $s \mapsto u(s)$. The covariant derivative is then:
$$
\nabla_{\dot{u}}Y := \dot{Y} - \frac{1}{4}g_{u}(\nabla Y, \nabla \dot{u}),
$$
where $\dot{u} = \frac{du}{dt}$. We refer to \cite[Section 7.5.1]{boyer_sasakian_2007} and \cite[Proposition 2.8]{franzinetti2020curvature} for more details.

By \cite[Lemma 5.1]{lempert2017riemannian}, the isometry $F:\Hspace(M_1, \mathcal{S}_1) \longrightarrow \Hspace(M_2, \mathcal{S}_2)$ preserves this expression, i.e.:
\begin{align*}
F_{\star\vert u} &\left(\dot{Y}(s) - \frac{1}{4} g_{1, u}(\nabla \dot{u} , \nabla Y )\right)\\
&= \frac{d}{dt}\left[F_{\star\vert u}(Y)\right] - \frac{1}{4} g_{2,F(u)} \left( \nabla (F_{\star \vert u}(\dot{u})), \nabla (F_{\star \vert u}(Y)) \right).
\end{align*}
Here, the metric $g_{1,u}$ is the Sasaki metric associated to $\mathcal{S}_{1,u}$ and similarly for $g_{2,F(u)}$. The gradients on the left hand side (resp. right hand side) are taken with respect to $g_{1,u}$ (resp. $g_{2,F(u)}$). For any two arbitrary functions $\psi_1, \psi_2 \in \Cinf_B(M_1)$ and choose $u$ such that $\dot{u}(0) = \psi_1 \in \Cinf_B(M_1)$ and $Y$ constantly equal to $\psi_2 \in \Cinf_B(M_1)$. Since $F_\star$ is independent on $u$, one has:
\begin{align}
\label{eq:LienMetriqueIsometrieSasakiDeriveeCovariante}
F_\star \left( g_{1,u} (\nabla \psi_1, \nabla \psi_2) \right) = g_{2, F(u)} \left( \nabla F_\star( \psi_1) , \nabla F_\star (\psi_2) \right).
\end{align}

\bigskip

Recall that $F_\star$ is a $\R$-linear map, thus choosing $\psi_1 = 1$ in \eqref{eq:LienMetriqueIsometrieSasakiDeriveeCovariante} reads:
$$
0 = g_{2, F(u)} \left( \nabla F_\star(1) , \nabla F_\star (\psi_2) \right), \quad \forall \psi_2 \in \Cinf_B(M_1),
$$
which in turn implies that $F_\star (1) \in \Cinf_B(M_1)$ must be constant. Also, integrating \eqref{eq:LienMetriqueIsometrieSasakiDeriveeCovariante} gives:
\begin{align}
\frac{1}{V_2} \int_{M_2} F_\star& \left( g_{1,u} (\nabla \psi_1, \nabla \psi_2) \right) \eta_2 \wedge d\eta_{2, F(u)}^n 
\\
&= \frac{1}{V_2} \int_{M_2} g_{2, F(u)} \left( \nabla F_\star( \psi_1) , \nabla F_\star (\psi_2) \right) \eta_2 \wedge d\eta_{2, F(u)}^n.
\label{eq:proofAffineIsometries1}
\end{align}
Here, $V_i$ stands for the volume $\Vol(M_i,\mathcal{S}_i)$, $i=1,2$. Since $F$ is an isometry, the left hand side above is:
\begin{equation}
\label{eq:proofAffineIsometries2}
\scal{1}{F_\star \left( g_{1,u} (\nabla \psi_1, \nabla \psi_2) \right)}_{F(u)}^{\Hspace(M_2,\mathcal{S}_2)} 
= \frac{1}{F_\star(1)}\scal{1}{g_{1,u} (\nabla \psi_1, \nabla \psi_2)}_u^{\Hspace(M_1,\mathcal{S}_1)},
\end{equation}
Now combining \eqref{eq:proofAffineIsometries1} and \eqref{eq:proofAffineIsometries2} then using integration by parts give:

\begin{align*}
\frac{1}{V_1 F_\star(1)} \int_{M_1} \left(\Delta_u \psi_1\right) \psi_2 \; \eta_1 \wedge d\eta_{1, u}^n 
&= \frac{1}{V_2} \int_{M_2} \left(\Delta_{F(u)} F_\star (\psi_1) \right) F_\star (\psi_2) \eta_2 \wedge d\eta_{2, F(u)}^n \\
&= \frac{1}{V_1} \int_{M_1} F_\star^{-1} \left(\Delta_{F(u)} F_\star (\psi_1) \right) \psi_2 \eta_1 \wedge d\eta_{1, u}^n.
\end{align*}
The last inequality above follows from the isometry condition with the inverse map $F^{-1}$. Here, $\Delta_u$ (resp. $\Delta_{F(u)}$) stands for the Riemannian Laplacian with respect to $g_{1,u}$ (resp. $g_{2,F(u)}$). By Riemannian Laplacian, we mean: $\Delta := - \star d \star d$, where $\star$ is the Hodge star operator.

This being true for any $\psi_2 \in \Cinf_B(M_1)$, we finally obtain:
\begin{align}
\label{Eq:AffineIsometryGivesCommutationWithLaplacian}
\frac{1}{F_\star (1)} F_\star\left( \Delta_u \psi \right)= \Delta_{F(u)} \left( F_\star (\psi)\right), \quad \forall \psi \in \Cinf_B(M_1).
\end{align}

The isometry condition for $F$ gives $F_\star (1) = \pm 1$. Since the Laplacians are non-negative operators, we must have: $F_\star(1) = 1$. This is reminiscent of Theorem \ref{Thm:INTROIsometriesKahler} more precisely, recall that in the Kähler setting, affine isometries are the ones with $b=0$, and thus we also have $F_\star(1) = 1$ in the Kähler setting.

On a Sasaki manifold $(M,\mathcal{S})$, we can consider different Laplacians. There is not only the Riemannian Laplacian but also the so called \emph{basic Laplacian} defined on the space of basic functions, $\Delta_B := - \overline{\star} d_B \overline{\star} d_B$ \cite[(7.2.4)]{boyer_sasakian_2007}, where $\overline{\star}(\cdot) := \star(\eta \wedge \cdot) = (-1)^r \iota_{\xi} (\star \cdot)$ on $r$-forms \cite[(7.2.2) and (7.2.3)]{boyer_sasakian_2007} and $d_B$ is the restriction of the usual differential to the space of basic functions \cite[Section 2.4]{boyer_sasakian_2007}. Note that if $\alpha$ is a basic form, then $\overline{\star} \alpha$ is so.
%Cf page 289 Nov2019

These two Laplacians coincide on $\Cinf_B(M)$:

\begin{Prop}
\label{Prop:SasakianLaplacianBasicCoincide}
Given any $f\in \Cinf_B(M)$, we have: $\Delta f = \Delta_B f$.
\end{Prop}

A proof of this fact can be found in \cite[Section 2.5]{nitta2009diameter} but we give here an alternative proof.
\begin{proof}
Fix $f\in \Cinf_B(M)$,
\begin{align*}
\Delta_B f 
&= - \overline{\star} d_B \overline{\star} d_B f 
= - \overline{\star} d_B \left( -\iota_\xi \star df \right)  
= - \overline{\star}\left( \iota_\xi d \star f - \mathcal{L}_\xi \star df \right) \\
&= -\overline{\star} \left(\iota_\xi \star \star d \star df  - \Lie_\xi \star df \right) 
= \overline{\star}\iota_\xi \star \Delta f + \overline{\star} \Lie_\xi \star df
= \Delta f +\overline{\star} \Lie_\xi \star df.
\end{align*}
Here, $\mathcal{L}$ stands for the Lie derivative, we used Cartan formula : $\Lie = d \iota + \iota d$. We also used that since a Sasaki manifold is odd dimensional, we have $\star \star d \star df = d\star df$. In the last inequality, we used: $\overline{\star}\,\overline{\star} = (-1)^r$ over $r$-forms.
%Cf Nov2019 Page 287 : Indeed, let $\alpha$ be a $r$-form, $\star \star \alpha = (-1)^{r(2n+1-r)}\alpha$

Now the proof of the claim ends by noticing that the last term is $0$, indeed, it is equal to 
$$
\star \left( \eta \wedge \Lie_\xi \left( \star df \right) \right) \underbrace{=}_{\Lie_\xi \eta = 0} - \star \Lie_\xi \left(\eta\wedge \star df \right),
$$
since $df$ is basic so $\star df$ contains $\eta$ thus $\eta \wedge \star df = 0$. More precisely, choose a coordinate system $(t, x_1, \cdots x_n, y_1, \cdots, y_n)$ for which $\eta = dt$, $df$ lies in $\textit{span}(dx_1, \cdots dx_n, dy_1, \cdots, dy_n)$ thus $dt \wedge \star df = 0$.
\end{proof}

Proposition \ref{Prop:SasakianLaplacianBasicCoincide} together with \eqref{Eq:AffineIsometryGivesCommutationWithLaplacian} then give Theorem \ref{Th:INTRO_AffineIsometries}:

\begin{proof}[Proof of Theorem \ref{Th:INTRO_AffineIsometries}]
Let $f\in\Cinf_B(M_1)$ be an eigenfunction, i.e. $\Delta_{u,B} f = \lambda f$ for some $\lambda \in \R$. Then \eqref{Eq:AffineIsometryGivesCommutationWithLaplacian} implies that $F_\star(f)$ is an eigenfunction of $\Delta_{F(u),B}$ with respect to the eigenvalue $\frac{\lambda}{F_\star (1)}$.
\end{proof}

\begin{footnotesize}
\addcontentsline{toc}{section}{Bibliography}
\bibliographystyle{abbrv}
\bibliography{/home/thomas/Documents/LaTex/Bibliographie}
\end{footnotesize}
\end{document}